\documentclass[11pt,twoside,leqno]{article}

\usepackage{amsthm,amsmath,amsfonts,amsbsy}
\usepackage[numbers]{natbib}
\usepackage{hyperref}
\usepackage{fancyhdr}
\long\def\symbolfootnote[#1]#2{\begingroup%
\def\thefootnote{\fnsymbol{footnote}}\footnote[#1]{#2}\endgroup}

\makeatletter
\renewcommand\section{\@startsection {section}{1}{\z@}%
                                   {-3.5ex \@plus -1ex \@minus -.2ex}%default
                                   {2.3ex \@plus.2ex}%default
                                   {\centering\normalfont\bfseries}}
\makeatother

\newtheoremstyle{plaincaps}%
    {1em}%
    {1em}%
    {\it}%
    {\parindent}%
    {\sc}%
    {.}%
    {1ex}%
    {}%
\theoremstyle{plaincaps}
\newtheorem{theorem}{Theorem}[section]
\newtheorem{condition}[theorem]{Condition}
\newtheorem{corollary}[theorem]{Corollary}
\newtheorem{lemma}[theorem]{Lemma}
\newtheorem{proposition}[theorem]{Proposition}
\newtheorem{definition}[theorem]{Definition}

\newtheoremstyle{examplecaps}%
    {1ex}%
    {1ex}%
    {\rm}%
    {\parindent}%
    {\sc}%
    {.}%
    {1ex}%
    {}%
\theoremstyle{examplecaps}
\newtheorem{example}[theorem]{Example}

\numberwithin{equation}{section}

\newcommand{\eps}{\varepsilon}

\newcommand{\dto}{\rightsquigarrow}
\newcommand{\vto}{\stackrel{v}{\to}}

\newcommand{\floor}[1]{\lfloor#1\rfloor}

\newcommand{\RV}{\mathcal{R}}

\newcommand{\law}{\mathcal{L}}
\newcommand{\sign}{\mathop{\mathrm{sign}}}

\newcommand{\E}{\mathrm{E}}
\renewcommand{\Pr}{\mathrm{P}}

\newcommand{\rmd}{\mathrm{d}}

\newcommand{\du}{\mathrm{d}u}
\newcommand{\dv}{\mathrm{d}v}
\newcommand{\dx}{\mathrm{d}x}

\newcommand{\dz}{\mathrm{d}z}

\newcommand{\EE}{\mathbb{E}}
\newcommand{\RR}{\mathbb{R}}
\renewcommand{\SS}{\mathbb{S}}
\newcommand{\VV}{\mathbb{V}}
\newcommand{\WW}{\mathbb{W}}

\newcommand{\ZZ}{\mathbb{Z}}

\newcommand{\M}{\mathcal{M}}

\newcommand{\BFTC}{\mbox{\sc bftc}}

\newcommand{\1}{\boldsymbol{1}}

\begin{document}
\thispagestyle{plain}
\begin{center}
{\bf MULTIVARIATE REGULAR VARIATION OF HEAVY-TAILED MARKOV CHAINS\symbolfootnote[0]{This version: \today.}} \medskip

{\sc Johan Segers}\symbolfootnote[1]{Institut de statistique, Universit\'e catholique de Louvain, Voie du Roman Pays 20, B-1348 Louvain-la-Neuve, Belgium; E-mail: segers@stat.ucl.ac.be. Supported by the Netherlands Organization for Scientific Research (NWO) in the form of a VENI grant and by the IAP research network grant nr.\ P5/24 of the Belgian government (Belgian Science Policy).} \medskip

{\em Universit\'e catholique de Louvain}
\end{center}

\begin{quote}\footnotesize
{\bf Abstract.} 
The upper extremes of a Markov chain with regulary varying stationary marginal distribution are known to exhibit under general conditions a multiplicative random walk structure called the tail chain. More generally, if the Markov chain is allowed to switch from positive to negative extremes or vice versa, the distribution of the tail chain increment may depend on the sign of the tail chain on the previous step. But even then, the forward and backward tail chain mutually determine each other through a kind of adjoint relation. As a consequence, the finite-dimensional distributions of the Markov chain are multivariate regularly varying in a way determined by the back-and-forth tail chain. An application of the theory yields the asymptotic distribution of the past and the future of the solution to a stochastic difference equation conditionally on the present value being large in absolute value.


\end{quote}

% ==============================================================================
\section{Introduction}
\label{S:intro}

Consider a discrete-time, real-valued random process $\{ X_t : t = 0, 1, 2, \ldots \}$ defined by the recursive equation
\begin{equation}
\label{E:MC:1}
	X_t = \Psi ( X_{t-1}, \eps_t ), \qquad t = 1, 2, \ldots,
\end{equation}
where
\begin{equation}
\label{E:MC:2}
\mbox{\begin{minipage}[t]{0.90\textwidth}
\begin{itemize}
\item[(i)] $\eps_1, \eps_2, \ldots$ are independent and identically distributed random elements of a measurable space $(\mathbb{S}, \mathcal{S})$ and independent of $X_0$;
\item[(ii)] $\Psi$ is a measurable function from $\RR \times \SS$ to $\RR$.
\end{itemize}
\end{minipage}}
\end{equation}
If the process $\{ X_t \}$ happens to be stationary, it will be assumed to be defined for all integer $t$. The distribution of $X_0$ is assumed to be regularly varying.

Aim of the paper is to find the weak limits of the finite-dimensional distributions of the process conditionally on $X_0$ being large in absolute value. More precisely, we will investigate the weak limits, called the forward tail chain, of vectors of the form $(X_0, \ldots, X_t)$ given $|X_0|$ exceeds a high threshold. If in addition the process is stationary, we will extend to this to find the back-and-forth tail chain, which corresponds to the weak limits of vectors of the form $(X_{-s}, \ldots, X_t)$ given $|X_0|$ is large. As a consequence, the finite-dimensional distributions of the process $\{ X_t \}$ will under quite general conditions be found to be multivariate regularly varying.

The process $\{ X_t \}$ is obviously a discrete-time homogeneous Markov chain. Every homogeneous discrete-time Markov chain $\{ X_t \}$ can be represented as in \eqref{E:MC:1}--\eqref{E:MC:2}. Of course, for a given Markov chain $\{ X_t \}$ the above representation is not unique. Still, in examples, the way in which Markov chains are defined is often through a recursive equation; all examples in \citet[][pp.~126--127]{Goldie91}, for instance, are of this type. The chain is stationary if and only if $X_1 = \Psi(X_0, \eps_1)$ and $X_0$ are equal in law.

In \citet{Smith92} and \citet{Perfekt94}, excursions of a Markov chain over a high threshold are shown to behave asymptotically and under quite general conditions as a (multiplicative) random walk. The theory has been extended to multivariate Markov chains in \citet{Perfekt97} and to higher-order Markov chains in \citet{Yun98, Yun00}. The random-walk representation is useful from a statistical perspective because it gives a handle on how to model the extremes of certain time series \citep{BC00, CST97, STC97}. A useful, well-investigated class of processes for which the random walk structure is quite revealing are the stationary solutions to certain stochastic difference equations, including squared autoregressive conditionally heteroskedastic (ARCH) processes as a special case \citep{BDM02b, GdHP04, dHRRdV89}.

A limitation of the theory of \citet{Smith92} and \citet{Perfekt94} is that it excludes Markov chains for which extreme values can switch from the upper to the lower tail or vice versa, as can be observed for instance in time series of logreturns of prices of financial securities in periods of high volatility. For such Markov chains with tail switching potential, the random walk representation of excursions over high thresholds breaks down in the sense that the distribution of the multiplicative increment now depends in general on the sign of the chain on the previous step. In \citet{BC03}, a more general representation is postulated, involving in fact four transition mechanisms rather than one, corresponding to the four cases of transitions from and to upper or lower extreme states.

%In all cases, the basic ingredient in the theory of extremes of stationary Markov chains is the so-called forward tail chain, which corresponds to the asymptotic distribution of the process given that its initial value is large. Similarly, the backward tail chain corresponds to the asymptotic distribution of the process given that it ends with a large value. Of course, the backward tail chain is nothing more than the forward tail chain of the time-reversed process, which is a Markov chain as well. 

The novelty of the paper is two-fold: first, to derive an elementary condition on the recursive mechanism \eqref{E:MC:1} under which the random walk representation or the more general tail-switching representation holds, including a simple description of the latter; second, in the stationary case, to study the joint distribution of the forward and backward tail chain, coined the {\em back-and-forth tail chain}. Besides the assumption that the distribution of $X_0$ is regularly varying, the only condition is a relatively easy-to-check statement on the asymptotic behaviour of $\Psi(x, \, \cdot \,)$ for large $|x|$.

The outline of the paper is as follows. The forward tail chain of a possibly non-stationary Markov chain is studied in section~\ref{S:forward}. Section~\ref{S:adjoint} describes a kind of adjoint relation between bivariate distributions that serves to define a class of processes, coined back-and-forth tail chains, in section~\ref{S:BFTC}, which in turn characterize the joint forward and backward tail chains of certain stationary Markov chains in section~\ref{S:backforth}. Multivariate regular variation and maximal domains of attraction are studied in section~\ref{S:MRV}. Finally, section~\ref{S:examples} provides some examples and counterexamples to the theory, including an application to stationary solutions of stochastic difference equations.

To conclude this section, let us fix some notations. For real $x$ and $y$, put $x \vee y = \max(x, y)$ and $x \wedge y = \min(x, y)$. For $0 < \alpha < \infty$ and $x \in \RR$, denote $(x)_+^\alpha = (x \vee 0)^\alpha$. The law of a random vector $\boldsymbol{X}$ is denoted by $\law(\boldsymbol{X})$; weak convergence of probability measures is denoted by $\dto$. The probability measure degenerate at a point $x$ is denoted by $\delta_x$. The indicator of an event $A$ is denoted by $\1(A)$.

% ==============================================================================
\section{Forward tail chain}
\label{S:forward}

Let $X_0, X_1, X_2, \ldots$ be a homogeneous Markov chain as in \eqref{E:MC:1} and \eqref{E:MC:2}, not necessarily stationary. The focus of this section is on the asymptotic distributions of the finite-dimensional distributions of the process conditionally on $|X_0|$ being large (Theorem~\ref{T:forward}). As a side result, the joint survival function of $(X_0, \ldots, X_t)$ is found to be multivariate regularly varying in all $2^{t+1}$ corners of $\RR^{t+1}$ (Condition~\ref{C:min}). Two conditions are required: Condition~\ref{C:RV} on the tails of $X_0$, and Condition~\ref{C:phi} on the asymptotics of $\Psi(x, \, \cdot \,)$ for large $|x|$. 

\begin{condition}
\label{C:RV}
There exists $0 < \alpha < \infty$ such that
\begin{equation}
\label{E:RV}
	\lim_{x \to \infty} \frac{\Pr(|X_0| > xy)}{\Pr(|X_0| > x)} = y^{-\alpha},
	\qquad 0 < y < \infty. 
\end{equation}
Moreover, there exists $0 \leq p \leq 1$ such that
\begin{equation}
\label{E:p}
	\lim_{x \to \infty} \frac{\Pr(X_0 > x)}{\Pr(|X_0| > x)} = p.
\end{equation}
Denote $p(+1) = p$ and $p(-1) = 1-p$.
\end{condition}

\begin{condition}
\label{C:phi}
There exists a measurable subset $\VV$ of $\SS$ with $\Pr(\eps_1 \in \VV) = 1$ such that for every $\sigma \in \{-1, 1\}$ with $p(\sigma) > 0$ and every $v \in \VV$ the following limit exists: 
\begin{equation}
\label{E:phi:V}
	\lim_{x \to \sigma \infty} x^{-1} \Psi(x, v) = \phi(v, \sigma).
\end{equation}
Moreover, if $\Pr[\phi(\eps_1, \sigma) = 0] > 0$ for some $\sigma$, then also $\Pr(\eps_1 \in \WW) = 1$, where $\WW$ is a measurable subset of $\VV$ such that for all $v \in \WW$,
\begin{equation}
\label{E:phi:W}
	\sup_{|y| \leq x} |\Psi(y, v)| = O(x), \qquad x \to \infty.
\end{equation}
\end{condition}

\begin{theorem}
\label{T:forward}
Let $\{ X_t \}$ be given by \eqref{E:MC:1}--\eqref{E:MC:2}. If Conditions~\ref{C:RV} and \ref{C:phi} hold, then for every integer $t \geq 0$, as $x \to \infty$,
\begin{equation}
\label{E:forward:1}
	\law \biggl( \frac{|X_0|}{x}, \frac{X_0}{|X_0|}, \frac{X_1}{|X_0|}, \ldots, \frac{X_t}{|X_0|} \bigg| |X_0| > x \biggr)
	\dto \law(Y, M_0, M_1, \ldots, M_t)
\end{equation}
with 
\begin{equation}
\label{E:Mj}
	M_j = M_{j-1} \phi(\eps_j, \sign M_{j-1}), \qquad j = 1, 2, \ldots,
\end{equation}
where $\phi(\, \cdot \, , 0) = 0$ and
\begin{equation}
\label{E:forward:2}
\mbox{\begin{minipage}[t]{0.90\textwidth}
\begin{itemize}
\item[(i)] $Y, M_0, \eps_1, \eps_2, \ldots$ are independent with $\eps_t$ as in \eqref{E:MC:2}(i);
\item[(ii)] $\Pr(Y > y) = y^{-\alpha}$ for $y \geq 1$;
\item[(iii)] $\Pr(M_0 = \sigma) = p(\sigma)$ for $\sigma \in \{-1, 1\}$.
\end{itemize}
\end{minipage}}
\end{equation}
\end{theorem}

\begin{proof}
The argument is by induction on $t$. The case $t = 0$ is a straightforward consequence of Condition~\ref{C:RV}. So let $t$ be a positive integer and let $f : \RR^{t+2} \to \RR$ be bounded and continuous. We have to show that
\begin{equation}
\label{E:forward:10}
	\lim_{x \to \infty}
	\E \biggl[
	f \biggl( \frac{|X_0|}{x}, \frac{X_0}{|X_0|}, \ldots, \frac{X_t}{|X_0|} \biggr)
	\, \biggl| \, |X_0| > x
	\biggr]
	= \E [ f(Y, M_0, \ldots, M_t) ].
\end{equation}
By \eqref{E:MC:1}, if $X_0 \neq 0$,
\[
	\frac{X_t}{|X_0|}
	= \frac{\Psi(X_{t-1}, \eps_t)}{|X_0|}
	= \frac{\Psi(x \frac{|X_0|}{x} \frac{X_{t-1}}{|X_0|}, \eps_t)}
	{x \frac{|X_0|}{x}}.
\]
Hence,
\begin{eqnarray}
\label{E:forward:20}
	\lefteqn{
	\E \biggl[
	f \biggl( \frac{|X_0|}{x}, \frac{X_0}{|X_0|}, \ldots, \frac{X_t}{|X_0|} \biggr)
	\, \biggl| \, |X_0| > x
	\biggr]
	} \\
	&=&
	\E \biggl[
	g_x \biggl( \frac{|X_0|}{x}, \frac{X_0}{|X_0|}, \ldots, \frac{X_{t-1}}{|X_0|} \biggr)
	\, \biggl| \, |X_0| > x
	\biggr] 
	\nonumber
\end{eqnarray}
where
\begin{equation}
\label{E:forward:gx}
	g_x(y, x_0, \ldots, x_{t-1})
	= \E \biggl[ f \biggl( y, x_0, \ldots, x_{t-1}, \frac{\Psi(xyx_{t-1}, \eps_t)}{xy} \biggr) 
	\biggr].
\end{equation}
Define
\begin{equation}
\label{E:forward:g}
	g(y_, x_0, \ldots, x_{t-1})
	= \E [ f (y, x_0, \ldots, x_{t-1}, x_{t-1} \phi(\eps_t, \sign x_{t-1}) ) ].
\end{equation}
By \eqref{E:Mj},
\begin{equation}
\label{E:forward:40}
	\E [ f(Y, M_0, \ldots, M_t) ] = \E [ g(Y, M_0, \ldots, M_{t-1}) ].
\end{equation}
In view of the identities \eqref{E:forward:20} and \eqref{E:forward:40}, the limit relation in \eqref{E:forward:10} will follow if we can show that
\begin{equation}
\label{E:forward:50}
	\E \biggl[
	g_x \biggl( \frac{|X_0|}{x}, \frac{X_0}{|X_0|}, \ldots, \frac{X_{t-1}}{|X_0|} \biggr)
	\, \biggl| \, |X_0| > x
	\biggr]
	\to \E [ g(Y, M_0, \ldots, M_{t-1}) ]
\end{equation}
as $x \to \infty$. In turn, \eqref{E:forward:50} will follow from the induction hypothesis and an extension of the continuous mapping theorem \citep[][Theorem~18.11]{vdV98} provided
\begin{equation}
\label{E:forward:55}
	\lim_{x \to \infty} g_x( y(x), x_0(x), \ldots, x_{t-1}(x) ) = g( y, x_0, \ldots, x_{t-1} )
\end{equation}
whenever $y(x) \to y$ and $x_i(x) \to x_i$ as $x \to \infty$ with $(y, x_0, \ldots, x_{t-1})$ ranging over a set $E \subset \RR^{t+1}$ with $\Pr[(Y, M_0, \ldots, M_{t-1}) \in E] = 1$. From the definitions of $g_x$ and $g$ in \eqref{E:forward:gx} and \eqref{E:forward:g}, respectively, equation~\eqref{E:forward:55} is implied by
\begin{equation}
\label{E:forward:60}
	\lim_{z \to \infty} \frac{\Psi(z w(z), v)}{z} = w \phi(v, \sign w)
\end{equation}
whenever $\lim_{z \to \infty} w(z) = w$ and where $w$ and $v$ range over sets that receive probability one by the distributions of $M_{t-1}$ and $\eps_1$, respectively. But \eqref{E:forward:60} is ensured by Condition~\ref{C:phi}: If $w \neq 0$, then \eqref{E:forward:60} follows from \eqref{E:phi:V}. The case $w = 0$ arises only if $\Pr(M_{t-1} = 0) > 0$; in turn, the latter can only occur if $\Pr[\phi(\eps_1, \sigma) = 0] > 0$ for some $\sigma$, and then \eqref{E:forward:60} is guaranteed by \eqref{E:phi:W}.
\end{proof}

\begin{corollary}
\label{C:min}
Under the conditions of Theorem~\ref{T:forward}, for every integer $t \geq 0$ and all $x_0, \ldots, x_t \in \RR$,
\begin{equation}
\label{E:min}
	\lim_{x \to \infty} \frac{\Pr( x_0 X_0 > x, \ldots, x_t X_t > x)}{\Pr(|X_0| > x)}
	= \E [ (x_0 M_0)_+^\alpha \wedge \cdots \wedge (x_t M_t)_+^\alpha ].
\end{equation}
\end{corollary}

\begin{proof}
If $x_0 = 0$, then both sides of \eqref{E:min} are zero. So assume $x_0 \neq 0$. By Theorem~\ref{T:forward} and the continuous mapping theorem, as $x \to \infty$,
\[
	\law \biggl( \frac{X_0}{x}, \frac{X_1}{x}, \ldots, \frac{X_t}{x} \, \bigg| \, |X_0| > x \biggr) \dto \law(Y M_0, Y M_1, \ldots, Y M_t).
\]
Hence, as the distribution of $Y$ is continuous, the limit on the left-hand side of \eqref{E:min} is equal to
\begin{eqnarray*}
	\lefteqn{
	\lim_{x \to \infty} \frac{\Pr(|X_0| > x/|x_0|)}{\Pr(|X_0| > x)}
	\Pr(x_0 X_0 > x, \ldots, x_t X_t > x \mid |X_0| > x/|x_0|)
	} \\
	&=& |x_0|^\alpha \Pr ( x_0 Y M_0 > |x_0|, x_1 Y M_1 > |x_0|, \ldots, x_t Y M_t > |x_0| ).
\end{eqnarray*}
It remains to show that this expression is equal to the right-hand side of \eqref{E:min}. Since the variable $Y$ is independent of $(M_0, \ldots, M_t)$ and $Y^{-\alpha}$ is uniformly distributed on $(0, 1)$, the above expression is equal to
\[
	|x_0|^\alpha \int_0^1 \Pr[ (x_0 M_0)_+^\alpha > |x_0|^\alpha u, (x_1 M_1)_+^\alpha > |x_0|^\alpha u, \ldots, (x_t M_t)_+^\alpha > |x_0|^\alpha u ] \du.
\]
Change variables $v = |x_0|^\alpha u$ and use $|M_0| = 1$ to simplify this to
\[
	\int_0^\infty 
	\Pr[ (x_0 M_0)_+^\alpha > v, \ldots, (x_t M_t)_+^\alpha > v ] \dv 
	= \E [ (x_0 M_0)_+^\alpha \wedge \cdots \wedge (x_t M_t)_+^\alpha ],
\]
as stated.
\end{proof}

In \citet[][Theorem~2.1]{GdHP04}, equation~\eqref{E:min} is derived for stationary solutions to certain stochastic difference equations and positive $x_i$.

% ==============================================================================
\section{An adjoint relation between bivariate distributions}
\label{S:adjoint}

The purpose of this section is to prepare the ground for the definition of back-and-forth tail chains in section~\ref{S:BFTC}, which in turn will show up in Theorem~\ref{T:backforth}, the main result of the paper. For $0 \leq p \leq 1$ and $0 < \alpha < \infty$, consider the following set of probability measures on $\{-1, 1\} \times \RR$:
\begin{equation}
\label{E:Mpa}
\mbox{\begin{minipage}[t]{0.90\textwidth}
$\M_{p,\alpha}$ is the set of distributions $\mu$ of bivariate random vectors $(I, M)$ with the property that
\begin{itemize}
\item[(i)] $\Pr(I = 1) = p = 1 - \Pr(I = -1)$;
\item[(ii)] $\E [ (\sigma M)_+^\alpha ] \leq \Pr(I = \sigma)$ for $\sigma \in \{-1, 1\}$.
\end{itemize}
\end{minipage}}
\end{equation}
Two measures $\mu, \mu^* \in \M_{p, \alpha}$ are said to be {\em adjoint} if, whenever $\law(I, M) = \mu$ and $\law(I^*, M^*) = \mu^*$,
\begin{equation}
\label{E:adjoint:min}
	\E [ (x I)_+^\alpha \wedge (y M)_+^\alpha ]
	= \E [ (x M^*)_+^\alpha \wedge (y I^*)_+^\alpha ],
	\qquad x, y \in \RR.
\end{equation}
Clearly, the relation ``\dots\ is adjoint to \dots'' is symmetric.

\begin{proposition}
\label{P:adjoint}
For every $\mu \in \M_{p, \alpha}$ there exists a unique $\mu^* \in \M_{p,\alpha}$ that is adjoint to $\mu$. If $\law(I, M) = \mu$ and $\law(I^*, M^*) = \mu^*$, then for $M^*$-integrable functions $f$ and for $\sigma \in \{-1, 1\}$,
\begin{eqnarray}
\label{E:adjoint:f}
	\lefteqn{ \E [ \1(I^* = \sigma) f(M^*) ] } \\
	&=& \E [ f(I / |M|) (\sigma M)_+^\alpha ] + \{ \Pr(I = \sigma) - \E [ (\sigma M)_+^\alpha ] \} f(0).
	\nonumber
\end{eqnarray}
\end{proposition}

\begin{proof}
Let $\mu \in \M_{p, \alpha}$ and $\law(I, M) = \mu$. Define a probability measure $\mu^*$ on $\{-1, 1\} \times \RR$ by
\begin{eqnarray*}
	\mu^*(\{(\sigma, 0)\})
	&=& \Pr(I = \sigma) - \E [ (\sigma M)_+^\alpha ] \\
	\mu^*(\{\sigma\} \times E )
	&=& \E \biggl[ \1\biggl( \frac{I}{|M|} \in E \biggr) (\sigma M)_+^\alpha \biggr]
\end{eqnarray*}
for $\sigma \in \{-1, 1\}$ and Borel sets $E \subset \RR \setminus \{0\}$. 

Let $(I^*, M^*)$ be a random pair with law $\mu^*$. By definition, equation~\eqref{E:adjoint:f} holds for measurable indicator functions $f$. By linearity, \eqref{E:adjoint:f} then also holds for measurable step functions, and therefore, by monotone convergence, for measurable nonnegative functions. Finally, by linearity, \eqref{E:adjoint:f} must hold for arbitrary $M^*$-integrable functions.

The measure $\mu^*$ belongs to $\M_{p, \alpha}$ because, by \eqref{E:adjoint:f}, $\Pr(I^* = \sigma) = \Pr(I = \sigma)$ and $\E [ (\sigma M^*)_+^\alpha ] = \Pr(I = \sigma, M \neq 0)$ for $\sigma \in \{-1, 1\}$. Moreover, $\mu^*$ solves \eqref{E:adjoint:min}, since the latter equation is a special case of \eqref{E:adjoint:f} applied to $\sigma = \sign y$ and $f(m^*) = (x m^*)_+^\alpha \wedge |y|^\alpha$.

To see that the solution of \eqref{E:adjoint:min} is necessarily unique, argue as follows. For $x, y \in \RR$,
\begin{eqnarray*}
	\E [ (xM^*)_+^\alpha \wedge (yI^*)_+^\alpha ]
	&=& \E [ \1(I^* = \sign y) (xM^*)_+^\alpha \wedge |y|^\alpha ] \\
	&=& \E \biggl[ 	\1(I^* = \sign y) 
					\int_0^{|y|^\alpha} \1 (xM^* > z^{1/\alpha}) \dz \biggr] \\
	&=& \int_0^{|y|^\alpha} \Pr(I^* = \sign y, \, x M^* > z^{1/\alpha}) \dz.
\end{eqnarray*}
Therefore, knowledge of the function $(x, y) \mapsto \E [ (xM^*)_+^\alpha \wedge (yI^*)_+^\alpha ]$ implies knowledge of $\mu^*$ on $\{-1, 1\} \times (\RR \setminus \{ 0 \})$. Since also $\mu^*(\{ 1 \} \times \RR) = p$ and $\mu^*(\{ -1 \} \times \RR) = 1 - p$, equation~\eqref{E:adjoint:min} uniquely determines $\mu^*$. 
\end{proof}

Let $\mu = \law(I, M)$ and $\mu^* = \law(I^*, M^*)$ be adjoint in $\M_{p, \alpha}$. By \eqref{E:adjoint:f},
\[
	\Pr(I^* = \sigma, \sign M^* = \tau) = \E[ \1(I = \tau) (\sigma M)_+^\alpha ], 
	\qquad \sigma, \tau \in \{-1, 1\},
\]
and in particular $\Pr(M^* \neq 0) = \E [ |M|^\alpha ]$ and vice versa. For $\sigma \in \{-1, 1\}$ such that $\Pr(I = \sigma) > 0$,
\begin{eqnarray}
\label{E:adjoint:c}
	\lefteqn{
	\E [f(M^* / I^*) \mid I^* = \sigma]
	} \\
	&=& \frac{1}{\Pr(I = \sigma)} 
	\E [ f(I/M) (\sigma M)_+^\alpha ] +
	\biggl( 1 - \frac{\E [ (\sigma M)_+^\alpha ]}{\Pr(I = \sigma)} \biggr) f(0), \nonumber
\end{eqnarray}
and vice versa. As a consequence,
\begin{equation}
\label{E:adjoint:d}
	\E [ f(M^* / I^*) ] = \E [ f(I / M) |M|^\alpha ] + (1 - \E [ |M|^\alpha ]) f(0).
\end{equation}
If additionally $\E [ |M|^\alpha ] = 1$, then $\E [ (\sigma M)_+^\alpha ] = \Pr(I = \sigma)$ for $\sigma \in \{-1, 1\}$ by~\eqref{E:Mpa}(ii), making the second terms on the right-hand sides of \eqref{E:adjoint:c} and \eqref{E:adjoint:d} vanish.

\begin{example}
\label{Ex:adjoint:ARCH}
Let $\mu = \law(I, IZ)$ where $I$ and $Z$ are independent random variables with $\Pr(I = 1) = 1/2 = \Pr(I = -1)$ and $\E [ |Z|^\alpha ] \leq 1$. Then $\mu \in \M_{1/2, \alpha}$ and by \eqref{E:adjoint:c}, its adjoint, $\mu^*$, is the distribution of $(I^*, I^* Z^*)$ where $I^*$ and $Z^*$ are independent random variables with $I^*$ equal in law to $I$ and with the distribution of $Z^*$ and $Z$ mutually determining each other via
\begin{equation}
\label{E:adjoint:ARCH}
	\E [ f(Z^*) ] = \E [ f(1/Z) |Z|^\alpha ] + (1 - \E [ |Z|^\alpha ]) f(0)
\end{equation}
and vice versa; see also Example~\ref{Ex:ARCH}.
\end{example}

\begin{example}
\label{Ex:adjoint:1}
A useful special case of Proposition~\ref{P:adjoint} arises if $p = 1$. If $\mu \in \M_{1, \alpha}$, then $\mu(\{1\} \times [0, \infty)) = 1$, so we can write $\mu = \delta_1 \otimes \nu$ for some probability measure $\nu$ on $[0, \infty)$ such that $\int_{[0, \infty)} x^\alpha \nu(\dx) \leq 1$. The adjoint of $\mu$ in $\M_{1, \alpha}$ is $\mu^* = \delta_1 \otimes \nu^*$ where $\nu^*$ is the probability measure on $[0, \infty)$ related to $\nu$ via
\begin{equation}
\label{E:adjoint:+:f}
	\E [ f(Z^*) ] = \E [ f(1/Z) Z^\alpha ] + (1 - \E[Z^\alpha]) f(0)
\end{equation}
and vice versa, where $\law(Z) = \nu$ and $\law(Z^*) = \nu^*$. 
Some examples of pairs $(\nu, \nu^*)$ are the following:
\begin{itemize}
\item For every $\alpha$, if $\nu$ is concentrated on $\{0, 1\}$, then $\nu^* = \nu$.
\item If $\alpha = 1$ and $\nu$ is the distribution of a unit exponential random variable, then $\nu^*$ is the distribution of the reciprocal of the sum of two independent unit exponential random variables.
\item If $\alpha = 1$ and $\nu$ is the distribution of a lognormal random variable with unit expectation, then $\nu^* = \nu$.
\end{itemize}
Except for a change of sign, the case $p = 0$ is similar to the case $p = 1$.
\end{example}

\begin{example}
\label{Ex:adjoint:Smith}
Let $\mu$ be the law of $(I, M)$ where $\Pr(I = 1) = p = 1 - \Pr(I = -1)$ for some $0 < p < 1$ and
\[
	\law( M/I \mid I = \pm 1) = p_1 \delta_1 + p_0 \delta_0 + p_{-1} \delta_{-\{(1-p)/p\}^{\pm 1/\alpha}}
\]
for some $0 < \alpha < \infty$ and some $p_{-1}, p_0, p_1 \in [0, 1]$ with $p_{-1} + p_0 + p_1 = 1$. Then $\mu \in \M_{p, \alpha}$ and the adjoint of $\mu$ in $\M_{p, \alpha}$ is $\mu$ itself; see also Example~\ref{Ex:Smith}.
\end{example}

% ==============================================================================
\section{Back-and-forth tail chains}
\label{S:BFTC}

In this section, the preparation of Theorem~\ref{T:backforth} is continued through the study of a certain class of discrete-time processes. Recall the set $\M_{p, \alpha}$ in \eqref{E:Mpa} and the adjoint relation in \eqref{E:adjoint:min}. 

\begin{definition}
\label{D:BFTC}
A discrete-time process $\{ M_t : t \in \ZZ \}$ is said to be a back-and-forth tail chain with index $0 < \alpha < \infty$ and forward transition law $\mu \in \M_{p, \alpha}$, notation $\BFTC(\alpha, \mu)$, if
\begin{itemize}
\item[(i)] $\law(M_0, M_1) = \mu$; 
\item[(ii)] $\law(M_0, M_{-1}) = \mu^*$ is adjoint to $\mu$ in $\M_{p, \alpha}$;
\item[(ii)] for all integer $t \geq 1$ and all real $x_{t-1}, x_{t-2}, \ldots$, 
\begin{eqnarray*}
	\lefteqn{
	\law(M_t \mid M_{t-1} = x_{t-1}, M_{t-2} = x_{t-2}, \ldots)
	} \\
	&=& \left\lbrace \begin{array}{l@{\quad}l}
	\law(x_{t-1} M_1/M_0 \mid M_0 = \sign x_{t-1}) & \mbox{if $x_{t-1} \neq 0$,} \\
	\delta_0 & \mbox{if $x_{t-1} = 0$;}
	\end{array} \right.
\end{eqnarray*}
\item[(iii)] for all integer $t \geq 1$ and all real $x_{-t+1}, x_{-t+2}, \ldots$,
\begin{eqnarray*}
	\lefteqn{
	\law(M_{-t} \mid M_{-t+1} = x_{-t+1}, M_{-t+2} = x_{-t+2}, \ldots)
	} \\
	&=& \left\lbrace \begin{array}{l@{\quad}l}
	\law(x_{-t+1} M_{-1}/M_0 \mid M_0 = \sign x_{-t+1}) & \mbox{if $x_{-t+1} \neq 0$,} \\
	\delta_0 & \mbox{if $x_{-t+1} = 0$.}
	\end{array} \right.
\end{eqnarray*}
\end{itemize}
\end{definition}

If $\mu = \delta_1 \otimes \nu \in \M_{1, \alpha}$ as in Example~\ref{Ex:adjoint:1}, we also write $\BFTC(\alpha, \nu)$.

Clearly, $\{ M_t : t \in \ZZ \}$ is a $\BFTC(\alpha,\mu)$ if and only if $\{ M_{-t} : t \in \ZZ \}$ is a $\BFTC(\alpha,\mu^*)$. Necessarily $\Pr(M_0 = 1) = p = 1 - \Pr(M_0 = -1)$. 

Let $0 < p < 1$. For $\mu \in \M_{p, \alpha}$ with adjoint $\mu^* \in \M_{p, \alpha}$, a process $\{ M_t \}$ is a $\BFTC(\alpha,\mu)$ if and only if it admits the following distributional representation: for integer $t \geq 1$, the variables $M_t$ and $M_{-t}$ are recursively given by
\begin{eqnarray}
\label{E:BFTC:AB:1}
	M_t &=&
	\left\lbrace \begin{array}{l@{\quad}l}
		M_{t-1} A_t & \mbox{if $M_{t-1} > 0$,} \\
		0 & \mbox{if $M_t = 0$,} \\
		M_{t-1} B_t & \mbox{if $M_{t-1} < 0$;}
	\end{array} \right. \\
	M_{-t} &=&
	\left\lbrace \begin{array}{l@{\quad}l}
		M_{-t+1} A_{-t} & \mbox{if $M_{-t+1} > 0$,} \\
		0 & \mbox{if $M_{-t+1} = 0$,} \\
		M_{-t+1} B_{-t} & \mbox{if $M_{-t+1} < 0$;}
	\end{array} \right. \nonumber
\end{eqnarray}
here
\begin{equation}
\label{E:BFTC:AB:2}
\mbox{\begin{minipage}[t]{0.9\textwidth}
\begin{itemize}
\item[(i)] $M_0, A_1, A_{-1}, A_2, A_{-2}, \ldots, B_1, B_{-1}, B_2, B_{-2}, \ldots$ are independent;
\item[(ii)] $\Pr(M_0 = 1) = p = 1 - \Pr(M_0 = -1)$;
\item[(iii)] $\law(A_t) = \law(A_1)$, $\law(A_{-t}) = \law(A_{-1})$, $\law(B_t) = \law(B_1)$ and $\law(B_{-t}) = \law(B_{-1})$ for integer $t \geq 1$;
\item[(iv)] $\law(M_0, M_1) = \mu$ and $\law(M_0, M_{-1}) = \mu^*$.
\end{itemize}
\end{minipage}}
\end{equation}
The laws $A_{-1}$ and $B_{-1}$ in \eqref{E:BFTC:AB:2} can be expressed in terms of $(\alpha,\mu)$ through \eqref{E:adjoint:c}: since $\law(A_{-1}) = \law(M_{-1} \mid M_0 = 1)$ and $\law(B_{-1}) = \law(-M_{-1} \mid M_0 = -1)$, for suitably integrable $f$,
\begin{equation}
\label{E:BFTC:AB:4}
	\begin{array}{rcl}
	\E [ f(A_{-1}) ]
	&=& \displaystyle \frac{1}{p} \E \biggl[ f \biggl( \frac{M_0}{M_1} \biggr) (M_1)_+^\alpha \biggr]
	+ \biggl( 1 - \frac{\E[(M_1)_+^\alpha]}{p} \biggr) f(0), \\[1ex]
	\E [ f(B_{-1}) ]
	&=& \displaystyle \frac{1}{1-p} \E \biggl[ f \biggl( \frac{M_0}{M_1} \biggr) (-M_1)_+^\alpha \biggr]
	+ \biggl( 1 - \frac{\E[(-M_1)_+^\alpha]}{1-p} \biggr) f(0).
	\end{array}
\end{equation}

On the other hand, if $p = 1$ and $\mu = \delta_1 \otimes \nu \in \M_{1, \alpha}$ with adjoint $\mu^* = \delta_1 \otimes \nu^*$, a process $\{ M_t \}$ is a $\BFTC(\alpha, \nu)$ if and only if it admits the following distributional representation:
\begin{equation}
\label{E:BFTC:A:1}
	\begin{array}[t]{rcl}
	M_0 &=& 1, \\[1ex]
	M_{\pm t} &=& \prod_{i=1}^t A_{\pm i}, \qquad \mbox{integer $t \geq 1$},
	\end{array}
\end{equation}
where
\begin{equation}
\label{E:BFTC:A:2}
\mbox{\begin{minipage}[t]{0.90\textwidth}
\begin{itemize}
\item[(i)] $A_1, A_{-1}, A_2, A_{-2}, \ldots$ are independent;
\item[(ii)] $\law(A_t) = \nu$ and $\law(A_{-t}) = \nu^*$ for integer $t \geq 1$.
\end{itemize}
\end{minipage}}
\end{equation}
The case $p = 0$ is similar to the case $p = 0$.

The above representations imply that if $\{ M_t \}$ is a $\BFTC(\alpha, \mu)$, then
\begin{equation}
\label{E:BFTC:moment}
	\E [ |M_t|^\alpha ] \leq 1, \qquad t \in \ZZ.
\end{equation}

The following proposition states a remarkable identiy for back-and-forth tail chains. This identity will play a crucial role in the proof of Theorem~\ref{T:backforth}.

\begin{proposition}
\label{P:BFTC}
Let $\mu \in \M_{p, \alpha}$ and let $\{ M_t : t \in \ZZ\}$ be a $\BFTC(\alpha,\mu)$. For all integers $s \geq 1$ and $t \geq 0$ and for every bounded, measurable function $f : \RR^{s+t+1} \to \RR$ such that $f(x_{-s}, \ldots, x_t) = 0$ as soon as $x_{-s} = 0$, the numbers
\begin{equation}
\label{E:BFTC}
	\E	\biggl[
			f \biggl( \frac{M_{-s+i}}{|M_i|}, \ldots, \frac{M_{t+i}}{|M_i|} \biggr)
			|M_i|^\alpha
		\biggr],
	\qquad i = 0, \ldots, s
\end{equation}
are all the same.
\end{proposition}

\begin{proof}
Note that it is sufficient to show that the numbers in \eqref{E:BFTC} corresponding to $i = 0$ and $i = 1$ are the same, that is,
\begin{equation}
\label{E:BFTC:1}
	\E [ f(M_{-s}, \ldots, M_t) ]
	= \E \biggl[
		f \biggl( \frac{M_{-s+1}}{|M_1|}, \ldots, \frac{M_{t+1}}{|M_1|} \biggr)
		|M_1|^\alpha
	\biggr].
\end{equation}
For, if \eqref{E:BFTC:1} is true for all integer $s \geq 1$ and $t \geq 0$, then an application of \eqref{E:BFTC:1} with $s$ and $t$ replaced by respectively $s-j$ and $t+j$ for some $j = 0, \ldots, s-1$ shows that the expectation in \eqref{E:BFTC} corresponding to $i = j$ is equal to the expectation in \eqref{E:BFTC} corresponding to $i = j + 1$.

We focus on the case $0 < p < 1$, the proofs for the cases $p = 0$ or $p = 1$ being similar but simpler. Without loss of generality, we assume that the process $\{ M_t \}$ is given as in \eqref{E:BFTC:AB:1} and \eqref{E:BFTC:AB:2}. The proof of \eqref{E:BFTC:1} proceeds in two steps. \smallskip

{\it Step~1: $s = 1$.}
Let $f : \RR^{t+2} \to \RR$ be bounded and measurable and such that $f(0, x_0, \ldots, x_t) = 0$. We have to show that
\begin{equation}
\label{E:BFTC:2}
	\E [ f(M_{-1}, \ldots, M_t) ]
	= \E \biggl[ f \biggl( \frac{M_0}{|M_1|}, \ldots, \frac{M_{t+1}}{|M_1|} \biggr) |M_1|^\alpha \biggr].
\end{equation}
The proof is by induction on $t$. For $t = 0$, equation~\eqref{E:BFTC:2} reduces to \eqref{E:adjoint:f}. Let $t \geq 1$. Since $M_t = M_{t-1} A_t \1 (M_{t-1} > 0) + M_{t-1} B_t \1 (M_{t-1} < 0)$,
\begin{eqnarray*}
	\E [ f(M_{-1}, \ldots, M_t) ]
	&=& \E [ f(M_{-1}, \ldots, M_{t-1}, M_{t-1} A_t) \1 (M_{t-1} > 0) ] \\
	&& \mbox{} + \E [ f(M_{-1}, \ldots, M_{t-1}, M_{t-1} B_t) \1 (M_{t-1} < 0) ] \\
	&& \mbox{} + \E [ f(M_{-1}, \ldots, M_{t-1}, 0) \1 (M_{t-1} = 0) ].
\end{eqnarray*}
The variables $A_t$ and $B_t$ are independent of the vector $(M_{-1}, \ldots, M_{t-1})$. Hence, the above expectations do not change if $A_t$ and $B_t$ are replaced by $A_{t+1}$ and $B_{t+1}$, respectively. Condition on $A_{t+1}$ and $B_{t+1}$ and apply the induction hypothesis to see that the above expression is equal to
\begin{eqnarray*}
	&& \E \biggl[ 
	f \biggl( \frac{M_0}{|M_1|}, \ldots, \frac{M_t}{|M_1|}, A_{t+1} \frac{M_t}{|M_1|}
	\biggr) 
	\1 \biggl( \frac{M_t}{|M_1|} > 0 \biggr) |M_1|^\alpha 
	\biggr] \\
	&+& \E \biggl[ 
	f \biggl( \frac{M_0}{|M_1|}, \ldots, \frac{M_t}{|M_1|}, B_{t+1} \frac{M_t}{|M_1|}
	\biggr) 
	\1 \biggl( \frac{M_t}{|M_1|} < 0 \biggr) |M_1|^\alpha 
	\biggr] \\
	&+& \E \biggl[ 
	f \biggl( \frac{M_0}{|M_1|}, \ldots, \frac{M_t}{|M_1|}, 0 \biggr) 
	\1 \biggl( \frac{M_t}{|M_1|} = 0 \biggr) |M_1|^\alpha. 
	\biggr]
\end{eqnarray*}
As $M_{t+1} = M_t A_{t+1} \1 (M_t > 0) + M_t B_{t+1} \1 (M_t < 0)$, the above expression can be simplified to the right-hand side of \eqref{E:BFTC:2}, as was to be shown. \smallskip

{\it Step~2: general $s$ and $t$.}
We proceed by induction on $s$. The case $s = 1$ was treated in step~1. Let $s \geq 2$. Since $M_{-s} = M_{-s+1} A_{-s} \1(M_{-s+1} > 0) + M_{-s+1} B_{-s} \1(M_{-s+1} < 0)$ and since $f(M_{-s}, \ldots, M_t) \1(M_{-s} = 0) = 0$,
\begin{eqnarray*}
%	\lefteqn{
	\E [ f(M_{-s}, \ldots, M_t) ]
%	} \\
	&=& \E [ f(M_{-s+1} A_{-s}, M_{-s+1}, \ldots, M_t) \1 (M_{-s+1} > 0) ] \\
	&& \mbox{} + 
	\E [ f(M_{-s+1} B_{-s}, M_{-s+1}, \ldots, M_t) \1 (M_{-s+1} < 0) ].
\end{eqnarray*}
Note that $(A_{-s}, B_{-s})$ is independent of $(M_{-s+1}, \ldots, M_t)$. Conditionally on $A_{-s}$ and $B_{-s}$, apply the induction hypothesis on the vector $(M_{-s+1}, \ldots, M_t)$ to rewrite the above expression as
\begin{eqnarray*}
	&& \E	\biggl[ 
			f	\biggl(
					\frac{M_{-s+2}}{|M_1|} A_{-s}, \frac{M_{-s+2}}{|M_1|}, \ldots,
					\frac{M_{t+1}}{|M_1|}
				\biggr)
			\1 \biggl( \frac{M_{-s+2}}{|M_1|} > 0 \biggr) |M_1|^\alpha
		\biggr] \\
	&+&
	\E	\biggl[ 
			f	\biggl(
					\frac{M_{-s+2}}{|M_1|} B_{-s}, \frac{M_{-s+2}}{|M_1|}, \ldots,
					\frac{M_{t+1}}{|M_1|}
				\biggr)
			\1 \biggl( \frac{M_{-s+2}}{|M_1|} < 0 \biggr) |M_1|^\alpha
		\biggr] \\
\end{eqnarray*}
The two expectations in the above display do not change if $A_{-s}$ and $B_{-s}$ are replaced by $A_{-s+1}$ and $B_{-s+1}$, respectively. By definition of $M_{-s+1}$, the above expression is then indeed equal to the right-hand side of \eqref{E:BFTC:1}.
\end{proof}

A special case of equation~\eqref{E:BFTC} is that, for every integer $s \geq 1$ and every $x_0, \ldots, x_s \in \RR$, the numbers
\[
	\E [ (x_0 M_i)_+^\alpha \wedge (x_1 M_{i-1})_+^\alpha \wedge \cdots \wedge (x_s M_{i-s})_+^\alpha ],
	\qquad i = 0, \ldots, s
\]
are all the same.

% ==============================================================================
\section{Back-and-forth tail chains of stationary Markov chains}
\label{S:backforth}

From now on, the process $\{ X_t \}$ in \eqref{E:MC:1} and \eqref{E:MC:2} is assumed to be strictly stationary. A necessary and sufficient condition for stationarity is that 
\begin{equation}
\label{E:MC:3}
	\law(\Psi(X_0, \eps_1)) = \law(X_0).
\end{equation}
It may be highly non-trivial to find the law for $X_0$ that solves \eqref{E:MC:3}. But even when the stationary distribution does not admit an explicit expression, its tails may in many cases be found by the theory developed in \citet{Goldie91}.

If the process $\{ X_t \}$ is stationary, then by Kolmogorov's extension theorem, the range of $t$ can without loss of generality be assumed to be the set of all integers, $\ZZ$. Our aim is to extend Theorem~\ref{T:forward} and find the asymptotic distribution of $(X_{-s}, \ldots, X_t)$ conditionally on $|X_0| > x$ for all integer $s$ and $t$ (Theorem~\ref{T:backforth}). Recall the process $M_0, M_1, M_2, \ldots$ in \eqref{E:Mj} and \eqref{E:forward:2} and recall the set $\M_{p, \alpha}$ in \eqref{E:Mpa}.

\begin{lemma}
\label{L:M1}
Let $\{ X_t : t \in \ZZ\}$ be a stationary Markov chain with distribution determined by \eqref{E:MC:1}, \eqref{E:MC:2} and \eqref{E:MC:3}. If Conditions~\ref{C:RV} and \ref{C:phi} hold, then
\[
	\E [ (\sigma M_1)_+^\alpha ] \leq p(\sigma), \qquad \sigma \in \{-1, 1\},
\]
that is, $\law(M_0, M_1) \in \M_{p, \alpha}$.
\end{lemma}

\begin{proof}
Fix $\sigma \in \{-1, 1\}$. By \eqref{E:p},
\[
	\lim_{x \to \infty} \Pr( \sigma X_0 > x \mid |X_0| > x) = p(\sigma).
\]
Fix $\delta > 0$. By stationarity, equation~\eqref{E:RV}, and the case $t = 1$ in Theorem~\ref{T:forward}, 
\begin{eqnarray*}
	p(\sigma)
	&\geq& \lim_{x \to \infty} \Pr( |X_0| > \delta x, \sigma X_1 > x \mid |X_1| > x) \\
	&=& \lim_{x \to \infty} \frac{\Pr( |X_0| > \delta x )}{\Pr(|X_1| > x)}
	\Pr(\sigma X_1 > x \mid |X_0| > \delta x) \\
	&=& \delta^{-\alpha} \Pr[ \sigma M_1 > ( \delta Y )^{-1} ].
\end{eqnarray*}
Since the distribution of the random variable $Y^{-\alpha}$ is uniform on $(0, 1)$, the right-hand side in the previous display is equal to
\[
	\delta^{-\alpha} 
	\int_0^1 \Pr[ (\sigma M_1)_+^\alpha > \delta^{-\alpha} u ]	\du
	= \int_0^{\delta^{-\alpha}} \Pr[ (\sigma M_1)_+^\alpha > v ] \dv.
\]
For every $\delta > 0$, the integral on the right-hand side of the previous display is bounded from above by $p(\sigma)$; hence the same must be true for its limit as $\delta \to \infty$, which is $\E [ (\sigma M_1)_+^\alpha ]$.
\end{proof}

Recall Definition~\ref{D:BFTC} of a back-and-forth tail chain. Apart from stationarity, the conditions on the chain $\{ X_t \}$ are identical to those in Theorem~\ref{T:forward}.

\begin{theorem}
\label{T:backforth}
Let $\{ X_t : t \in \ZZ\}$ be a stationary Markov chain with distribution determined by \eqref{E:MC:1}, \eqref{E:MC:2} and \eqref{E:MC:3}. If Conditions~\ref{C:RV} and \ref{C:phi} hold, then for all nonnegative integer $s$ and $t$, as $x \to \infty$,
\begin{equation}
\label{E:backforth}
	\law \biggl( \frac{|X_0|}{x}, \frac{X_{-s}}{|X_0|}, \ldots, \frac{X_t}{|X_0|}
	\,\bigg|\, |X_0| > x \biggr)
	\dto \law (Y, M_{-s}, \ldots, M_t)
\end{equation}
with
\begin{equation}
\label{E:backforth:bis}
\mbox{\begin{minipage}[t]{0.90\textwidth}
\begin{itemize}
\item[(i)] $Y$ is independent of $\{ M_t : t \in \ZZ \}$;
\item[(ii)] $\Pr(Y > y) = y^{-\alpha}$ for $y \geq 1$;
\item[(iii)] $\{ M_t : t \in \ZZ \}$ is a $\BFTC(\alpha, \mu)$ where $\mu = \law(M_0, M_1)$ with $(M_0, M_1)$ as in \eqref{E:Mj} and \eqref{E:forward:2}.
\end{itemize}
\end{minipage}}
\end{equation}
\end{theorem}

\begin{proof}
By Lemma~\ref{L:M1}, the law of $(M_0, M_1)$ belongs to $\M_{p, \alpha}$. Let $Y$ and $\{ M_t : t \in \ZZ \}$ be as in \eqref{E:backforth:bis}. Equation~\eqref{E:backforth} will follow if we can show that
\begin{equation}
\label{E:backforth:2}
	\lim_{x \to \infty}
	\E \biggl[ f 
		\biggl( 
			\frac{|X_0|}{x}, \frac{X_{-s}}{|X_0|}, \ldots, \frac{X_t}{|X_0|} 
		\biggr)
	\,\bigg|\, |X_0| > x \biggr]
	= \E [ f(Y, M_{-s}, \ldots, M_t) ]
\end{equation}
for every bounded and uniformly continuous function $f : \RR^{s+t+2} \to \RR$. For $s = 0$, equation~\eqref{E:backforth:2} follows from Theorem~\ref{T:forward} and the representations of a $\BFTC(\alpha, \mu)$ in \eqref{E:BFTC:AB:1} and \eqref{E:BFTC:AB:2} for $0 < p < 1$ and $\eqref{E:BFTC:A:1}$ and \eqref{E:BFTC:A:2} for $p = 1$. For $s \geq 1$, a function $f : \RR^{s+t+2} \to \RR$ can be decomposed as
\begin{eqnarray*}
	f(y, x_{-s}, \ldots, x_t)
	&=& \{ f(y, x_{-s}, \ldots, x_t) - f(y, 0, x_{-(s-1)}, \ldots, x_t) \} \\
	&& \mbox{} + f(y, 0, x_{-(s-1)}, \ldots, x_t).
\end{eqnarray*}
Hence, if we can show \eqref{E:backforth:2} for integer $s \geq 1$ and for functions $f$ that satisfy the additional constraint $f(y, x_{-s}, \ldots, x_t) = 0$ as soon as $x_{-s} = 0$, then by an induction argument, \eqref{E:backforth:2} must be true for general $f$.

Take integer $s \geq 1$ and $t \geq 0$ and let $f : \RR^{s+t+2} \to \RR$ be a bounded and uniformly continuous function such that $f(y, x_{-s}, \ldots, x_t) = 0$ as soon as $x_{-s} = 0$. Fix $\delta > 0$. We have
\begin{eqnarray}
\label{E:backforth:delta}
	\lefteqn{
		\E 	\biggl[ 
			f \biggl( 
			\frac{|X_0|}{x}, \frac{X_{-s}}{|X_0|}, \ldots, \frac{X_t}{|X_0|}
			\biggr) 
			\, \bigg| \, |X_0| > x 
			\biggr]
	} \\
	&=& \E 	\biggl[ 
			f \biggl( 
			\frac{|X_0|}{x}, \frac{X_{-s}}{|X_0|}, \ldots, \frac{X_t}{|X_0|}
			\biggr) 
			\1 (|X_{-s}| \leq \delta x)
			\, \bigg| \, |X_0| > x 
			\biggr] \nonumber \\
	&& \mbox{} +
		\E 	\biggl[ 
			f \biggl( 
			\frac{|X_0|}{x}, \frac{X_{-s}}{|X_0|}, \ldots, \frac{X_t}{|X_0|}
			\biggr) 
			\1 (|X_{-s}| > \delta x)
			\, \bigg| \, |X_0| > x 
			\biggr]. \nonumber
\end{eqnarray}
The first term on the right-hand side of \eqref{E:backforth:delta} is bounded by
\[
	\sup \{ |f(y, x_{-s}, \ldots, x_t)| : 
	y > 1, \boldsymbol{x} \in \RR^{s+t+1}, |x_{-s}| \leq \delta \},
\]
which converges to zero as $\delta \downarrow 0$ by the assumptions on $f$. By stationarity, the second term on the right-hand side of \eqref{E:backforth:delta} is equal to
\[
	\frac{\Pr( |X_0| > \delta x )}{\Pr( |X_0| > x )}
	\E \biggl[
			f \biggl( 
			\frac{|X_s|}{x}, \frac{X_0}{|X_s|}, \ldots, \frac{X_{s+t}}{|X_s|}
			\biggr) 
			\1 (|X_s| > x)
			\, \bigg| \, |X_0| > \delta x 		
	\biggr].
\]
As $x \to \infty$, this converges to
\begin{eqnarray*}
	\lefteqn{
	\delta^{-\alpha}
	\E \biggl[
		f \biggl( 
		\delta Y |M_s|, \frac{M_0}{|M_s|}, \ldots, \frac{M_{s+t}}{|M_s|} 
		\biggr)
		\1 \{ |M_s| > ( \delta Y )^{-1} \}
	\biggr]
	} \\
	&=& \delta^{-\alpha} \int_0^1
	\E \biggl[
		f \biggl(
			\delta u^{-1/\alpha} |M_s|, 
			\frac{M_0}{|M_s|}, \ldots, \frac{M_{s+t}}{|M_s|}
		\biggr)
		\1 ( |M_s|^\alpha > \delta^{-\alpha} u )
	\biggr]
	\du \\
	&=& \int_0^{\delta^{-\alpha}}
	\E \biggl[
		f \biggl(
			v^{-1/\alpha} |M_s|,
			\frac{M_0}{|M_s|}, \ldots, \frac{M_{s+t}}{|M_s|} 
		\biggr)
		\1 ( |M_s|^\alpha > v )
	\biggr]
	\dv.
\end{eqnarray*}
By \eqref{E:BFTC:moment} and dominated convergence, the latter integral converges as $\delta \downarrow 0$ to
\[
	\int_0^\infty
	\E \biggl[
		f \biggl( 
		v^{-1/\alpha} |M_s|, \frac{M_0}{|M_s|}, \ldots, \frac{M_{s+t}}{|M_s|} 
		\biggr)
		\1 (|M_s|^\alpha > v)
	\biggr] \dv.
\]
By the preceding arguments, the integral above must be equal to the limit on the left-hand side of \eqref{E:backforth:2}; so it suffices to show that the integral above is also equal to the right-hand side of \eqref{E:backforth:2}.
 
Change variables $y = v^{-1/\alpha} |M_s|$ to rewrite the above integral as
\[
	\int_1^\infty
	\E \biggl[
		f \biggl(
			y, \frac{M_0}{|M_s|}, \ldots, \frac{M_{s+t}}{|M_s|}
		\biggr)
		|M_s|^\alpha 
	\biggr] 
	\rmd(-y^{-\alpha}).
\]
Apply Proposition~\ref{P:BFTC} to simplify this expression to
\[
	\int_1^\infty
	\E [ f (y, M_{-s}, \ldots, M_t) ] \rmd(-y^{-\alpha})
	= \E [ f(Y, M_{-s}, \ldots, M_t) ],
\]
as required.
\end{proof}

If $0 < p < 1$, then the tail chain $\{ M_t \}$ in \eqref{E:backforth:bis}(iii) admits the representation in equations \eqref{E:BFTC:AB:1} to \eqref{E:BFTC:AB:4} with
\[
	\begin{array}{rclcl}
	\law(X_1 / X_0 \mid X_0 > x) &\dto& \law(A_1) &=& \law(\phi(\eps_1, 1)), \\[1ex]
	\law(X_1 / X_0 \mid X_0 < -x) &\dto& \law(B_1) &=& \law(\phi(\eps_1, -1)), \\[1ex]
	\law(X_{-1} / X_0 \mid X_0 > x) &\dto& \law(A_{-1}), \\[1ex]
	\law(X_{-1} / X_0 \mid X_0 < -x) &\dto& \law(B_{-1})
	\end{array}
\]
as $x \to \infty$. On the other hand, if $p = 1$, then the tail chain $\{ M_t \}$ in \eqref{E:backforth:bis}(iii) admits the representation in \eqref{E:BFTC:A:1} and \eqref{E:BFTC:A:2} with $\law(A_{\pm 1})$ as above.

% ==============================================================================
\section{Multivariate regular variation}
\label{S:MRV}

Under the conditions of Theorem~\ref{T:backforth}, the finite-dimensional distributions of the chain $\{ X_t \}$ are multivariate regularly varying (Theorem~\ref{T:MRV} and Corollary~\ref{C:simpler}). We also identify the extreme value attractors of random vectors of the form $(\sigma_0 X_0, \ldots, \sigma_t X_t)$, where $\sigma_i \in \{-1, 1\}$ for $i = 0, \ldots, t$ (Corollaries~\ref{C:max} and \ref{C:l}).

Recall that a positive, measurable function $V$ defined in a neighbourhood of infinity is {\em regularly varying} of index $\tau \in \RR$ if $\lim_{x \to \infty} V(xy) / V(x) = y^\tau$ for all positive $y$;  notation $V \in \RV_\tau$. The law of a random vector $(\xi_1,\ldots,\xi_d)$ is {\em multivariate regularly varying} if there exists $V \in \RV_{-\alpha}$ with $\alpha > 0$ and a non-trivial Radon measure $\mu$ on $\EE_d = [-\infty, \infty]^d \setminus \{ \boldsymbol{0} \}$ such that
\begin{equation}
\label{E:MRV:def}
	\frac{1}{V(x)} \Pr[ (x^{-1} \xi_1, \ldots, x^{-1} \xi_d) \in \cdot \, ]
	\vto \mu, \qquad x \to \infty,
\end{equation}
the arrow $\vto$ denoting {\em vague convergence} of measures \citep[][section~3.4]{Resnick87}. Equivalent definitions of multivariate regular variation can be found for instance in \citet[][p.~69]{Resnick86} or \citet{BDM02a}. The limit measure $\mu$ in \eqref{E:MRV:def} is necessarily homogeneous of order $-\alpha$. In particular, hyperplanes perpendicular to the coordinate axes do not receive any mass, that is, $\mu(\{ \boldsymbol{x} : x_i = z \}) = 0$ for all $i = 1, \ldots, d$ and $z \in [-\infty, \infty] \setminus 0$.

According to Theorem~\ref{T:MRV} below, the law of $(X_0, \ldots, X_t)$ is multivariate regularly varying for every integer $t \geq 0$. The limit measure $\mu_t$ can be expressed in terms of the tail chain $\{ M_j \}$. To this end, it is convenient to partition $\EE_{t+1}$ in one of the following two ways,
\begin{equation}
\label{E:partition}
	\EE_{t+1} 
	= \bigcup_{i=0}^t \mathbb{F}_{t,i}
	= \bigcup_{i=0}^t \mathbb{G}_{t,i}
\end{equation}
where, for $i = 0, \ldots, t$,
\begin{eqnarray*}
	\mathbb{F}_{t, i} &=& \{ (x_0, \ldots, x_t) : x_0 = \ldots = x_{i-1} = 0, x_i \neq 0 \}, \\
	\mathbb{G}_{t, i} &=& \{ (x_0, \ldots, x_t) : x_i \neq 0, x_{i+1} = \ldots = x_t = 0 \}.
\end{eqnarray*}

\begin{theorem}
\label{T:MRV}
Under the conditions of Theorem~\ref{T:backforth}, for every integer $t \geq 0$ and as $x \to \infty$,
\begin{equation}
\label{E:MRV}
	\frac{1}{\Pr(|X_0| > x)} \Pr[ (x^{-1} X_0, \ldots, x^{-1} X_t) \in \cdot \, ]
	\vto \mu_t
\end{equation}
in $\EE_{t+1}$, where $\mu_t$ is determined by
\begin{eqnarray*}
	\int_{\mathbb{F}_{t,0}} f \rmd \mu_t
	&=& \int_0^\infty \E [ f(zM_0, \ldots, zM_t) ] \rmd(-z^{-\alpha}), \\
	\int_{\mathbb{F}_{t,i}} f \rmd \mu_t
	&=& \int_0^\infty 
	\E [ f(0, \ldots, 0, zM_0, \ldots, z M_{t-i}) \1 (M_{-1} = 0) ] 
	\rmd(-z^{-\alpha})
\end{eqnarray*}
for $i = 1, \ldots, t$, and also by
\begin{eqnarray*}
	\int_{\mathbb{G}_{t,t}} f \rmd \mu_t
	&=& \int_0^\infty \E [ f(zM_{-t}, \ldots, zM_0) ] \rmd(-z^{-\alpha}), \\
	\int_{\mathbb{G}_{t,i}} f \rmd \mu_t 
	&=& \int_0^\infty 
	\E [ f(zM_{-i}, \ldots, z M_0, 0, \ldots, 0) \1 (M_1 = 0) ] 
	\rmd(-z^{-\alpha})
\end{eqnarray*}
for $i = 0, \ldots, t-1$.
\end{theorem}

\begin{proof}
Let $f : \EE_{t+1} \to \RR$ be continuous and with compact support. We have to show that
\begin{equation}
\label{E:MRV:f}
	\lim_{x \to \infty} \frac{1}{\Pr(|X_0| > x)} \E [ f(x^{-1} X_0, \ldots, x^{-1} X_t) ] = \int f \rmd \mu_t,
\end{equation}
where the latter integral can be computed according to one of the two partitions in \eqref{E:partition} and the corresponding formulas in the theorem. 

The assumption that the support of $f$ is compact means that there exists $\eps > 0$ such that $f(x_0, \ldots, x_t) = 0$ if $|x_i| \leq \eps$ for all $i = 0, \ldots, t$. Take $0 < \delta \leq \eps$. We have
\begin{eqnarray*}
	\lefteqn{
	\E [ f(x^{-1} X_0, \ldots, x^{-1} X_t) ]
	} \\
	&=& \E [ f(x^{-1} X_0, \ldots, x^{-1} X_t) \1 ( |X_0| > \delta x ) ] \\
	&& \mbox{} + \sum_{i=1}^t 
	\E [ f(x^{-1} X_0, \ldots, x^{-1} X_t) 
	\1 ({\textstyle\bigvee_{j=0}^{i-1}} |X_j| \leq \delta x, |X_i| > \delta x) ].
\end{eqnarray*}
On the one hand,
\begin{eqnarray*}
	\lefteqn{
	\lim_{x \to \infty}
	\frac{1}{\Pr(|X_0| > x)}
	\E [ f(x^{-1} X_0, \ldots, x^{-1} X_t) \1 ( |X_0| > \delta x ) ]
	} \\
	&=& \lim_{x \to \infty}
	\frac{\Pr(|X_0| > \delta x)}{\Pr(|X_0| > x)}
	\E [ f(x^{-1} X_0, \ldots, x^{-1} X_t) \mid |X_0| > \delta x ] \\
	&=& \delta^{-\alpha}
	\E [ f(\delta Y M_0, \ldots, \delta Y M_t) ] \\
	&=& \delta^{-\alpha}
	\int_1^\infty \E[ f(\delta y M_0, \ldots, \delta y M_t) ] 
	\rmd (-y^{-\alpha}) \\
	&=& \int_\delta^\infty
	\E [ f(z M_0, \ldots, z M_t) ] \rmd(-z^{-\alpha}).
\end{eqnarray*}
On the other hand, for $i = 1, \ldots, t$, by a similar argument,
\begin{eqnarray*}
	\lefteqn{
	\lim_{x \to \infty} \frac{1}{\Pr(|X_0| > x)}
	\E [ f(x^{-1} X_0, \ldots, x^{-1} X_t) 
	\1 ({\textstyle\bigvee_{j=0}^{i-1}} |X_j| \leq \delta x, |X_i| > \delta x) ]
	} \\
	&=& \int_\delta^\infty
	\E [ f(z M_{-i}, \ldots, z M_{t-i})
	\1 ( {\textstyle\bigvee_{j=-i}^{-1}} |M_j| \leq \delta/z) ]
	\rmd(-z^{-\alpha}).
\end{eqnarray*}
We obtain that for arbitrary $\delta \in (0, \eps]$, the limit on the left-hand side of \eqref{E:MRV:f} is equal to
\begin{eqnarray*}
	&& \int_\delta^\infty \E [ f(z M_0, \ldots, z M_t) \rmd(-z^{-\alpha}) \\
	&+& \sum_{i = 1}^{t} \int_\delta^\infty \E [ f(z M_{-i}, \ldots, z M_{t-i}) ]
	\1 ( {\textstyle\bigvee_{j=-i}^{-1}} |M_j| \leq \delta / z ) ] \rmd(-z^{-\alpha}).
\end{eqnarray*}
Take the limit as $\delta \downarrow 0$ and apply the dominated convergence theorem, which is justified by \eqref{E:BFTC:moment} and the fact that $f$ is continuous and has compact support, to see that the limit on the left-hand side of \eqref{E:MRV:f} is also equal to
\begin{eqnarray*}
	&& \int_0^\infty \E [f (z M_0, \ldots z M_t) ] \rmd(-z^{-\alpha}) \\
	&+& \sum_{i = 1}^{t} \int_0^\infty \E [ f(0, \ldots, 0, z M_0, \ldots, z M_i)
	\1(M_{-1} = 0) ] \rmd(-z^{-\alpha}).
\end{eqnarray*}
The last expression coincides with the stated expressions for $\int f \rmd \mu_t$ via the partition of $\EE_{t+1}$ into $\mathbb{F}_{t,0} \cup \cdots \cup \mathbb{F}_{t,t}$. The proof of the expressions for $\int f \rmd \mu_t$ via $\mathbb{G}_{t,i}$ is completely similar.
\end{proof}

\begin{corollary}
\label{C:simpler}
Under the conditions of Theorem~\ref{T:backforth}, if $\Pr(M_1 = 0) = 0$ or, equivalently, $\E [ |M_{-1}|]^\alpha = 1$, then for $\mu_t$-integrable $f$,
\begin{equation}
\label{E:simpler:a}
	\int f \rmd \mu_t
	= \int_0^\infty \E [ f(zM_{-t}, \ldots, zM_0) ] \rmd(-z^{-\alpha});
\end{equation}
and if $\Pr(M_{-1} = 0) = 0$ or, equivalently, $\E [|M_1|^\alpha] = 1$, then
\begin{equation}
\label{E:simpler:b}
	\int f \rmd \mu_t
	= \int_0^\infty \E [ f(zM_0, \ldots, zM_t) ] \rmd(-z^{-\alpha}).
\end{equation}
\end{corollary}

\begin{corollary}
\label{C:max}
Under the conditions of Theorem~\ref{T:backforth}, for every integer $t \geq 0$ and all $x_0, \ldots, x_t \in \RR$,
\[
	l_t (x_0, \ldots, x_t)
	= \lim_{x \to \infty}
	\frac{\Pr( x_0 X_0 > x \mbox{ or } \ldots \mbox{ or } x_t X_t > x)}{\Pr(|X_0| > x)}
\]
is equal to
\[
	\E [ (x_0 M_0)_+^\alpha \vee \cdots \vee (x_t M_t)_+^\alpha ] 
	+ \sum_{i = 1}^{t} 
	\E [ \{ (x_i M_0)_+^\alpha \vee \cdots \vee (x_t M_{t-i})_+^\alpha \} \1(M_{-1} = 0) ]
\]
and also to
\[
	\E [ (x_0 M_{-t})_+^\alpha \vee \cdots \vee (x_t M_0)_+^\alpha ]
	+ \sum_{i=0}^{t-1}
	\E [ \{ (x_0 M_{-i})_+^\alpha \vee \cdots \vee (x_i M_0)_+^\alpha \} \1(M_1 = 0) ].
\]
\end{corollary}

As in Corollary~\ref{C:simpler}, the expressions for $l_t$ in Corollary~\ref{C:max} simplify considerably if $\Pr(M_{-1} \neq 0) = \E[|M_1|^\alpha] = 1$ or $\Pr(M_1 \neq 0) = \E[|M_{-1}|^\alpha] = 1$.

\begin{corollary}
\label{C:l}
Under the conditions of Theorem~\ref{T:backforth}, if the sequence $\{a_n\}$ is such that $\lim_{n \to \infty} n \Pr(|X_0| > a_n) = 1$, then for integer $t \geq 0$, $\sigma_i \in \{-1, 1\}$ and $x_i \in (0, \infty)$ ($i = 0, \ldots, t$),
\[
	\lim_{n \to \infty}
	\{ \Pr(\sigma_0 X_0 \leq a_n x_0, \ldots, \sigma_t X_t \leq a_n x_t) \}^n
	= \exp \{ - l_t( \sigma_0 / x_0, \ldots, \sigma_t / x_t) \}.
\]
\end{corollary}

% ==============================================================================
\section{Examples}
\label{S:examples}

%Example~\ref{Ex:Smith} illustrates the basic theory and Examples~\ref{Ex:Smith:counterexample} and \ref{Ex:dependent} describe some surprising phenomena that can occur if the first or the second part of Condition~\ref{C:phi} are not fulfilled; the first two examples are inspired by Example~3 on p.~43 in \citet{Smith92}. In Example~\ref{Ex:ARCH}, the joint extremes of the stationary solution to a stochastic difference equation, including squared ARCH processes, are described using the theory of the previous sections. 

Examples~\ref{Ex:ARCH} to \ref{Ex:Smith} illustrate the theory to ARCH(1) processes, more general stochastic difference equations, and a process inspired by Example~3 in \citet[][p.~43]{Smith92}, respectively. Examples~\ref{Ex:Smith:counterexample} and \ref{Ex:dependent} are counterexamples illustrating some surprising phenomena that can occur if the first or the second part of Condition~\ref{C:phi} are not fulfilled.

\begin{example}
\label{Ex:ARCH}
Let $0 < \beta < \infty$ and $0 < \lambda < 2 \mathrm{e}^\gamma$ where $\gamma$ is Euler'c constant, and consider the ARCH(1) process
\begin{equation}
\label{E:ARCH}
	X_t = (\beta + \lambda X_{t-1}^2)^{1/2} Z_t, \qquad t = 1, 2, \ldots
\end{equation}
where $X_0, Z_1, Z_2, \ldots$ are independent random variables and the common distribution of the $Z_t$ is standard normal. By \citet[][Theorem~5]{Kesten73} or \citet[][Theorem~4.1]{Goldie91} to the squared series $X_t^2$, a stationary solution to \eqref{E:ARCH} exists, and this stationary distribution satisfies Condition~\ref{E:RV} with $p = 1/2$ and $\alpha$ equal to the unique positive solution to the equation $\E [ \lambda^{\alpha / 2} |Z|^\alpha] = 1$; see also \citet[][section~8.4.2]{EKM97}. Condition~\ref{C:phi} is easily verified with $\phi(z, \pm 1) = \pm \lambda^{1/2} z$. The back-and-forth tail chain is given by Theorem~\ref{T:backforth} and, by symmetry of the normal distribution and Example~\ref{E:adjoint:ARCH}, it can be represented by
\begin{eqnarray*}
	M_t &=& M_0 \lambda^{t/2} Z_1 \cdots Z_t, \\
	M_{-t} &=& M_0 \lambda^{-t/2} Z_1^* \cdots Z_t^*,
\end{eqnarray*}
for integer $t \geq 1$, where $M_0, Z_1, Z_1^*, Z_2, Z_2^*, \ldots$ are independent random variables with $\Pr(M_0 = 1) = 1/2 = \Pr(M_0 = -1)$ and $Z_t$ standard normal and with the common law of the variables $Z_t^*$ related to the standard normal distribution via
\[
	\E [ f(Z_1^*) ] = \lambda^{\alpha/2} \E [ f(1 / Z_1) |Z_1|^\alpha ],
\]
for $Z_1^*$-integrable functions $f$. By Corollary~\ref{C:max}, for all integer $t \geq 0$ and all real $x_0, \ldots, x_t$,
\begin{eqnarray*}
	\lefteqn{
	\lim_{x \to \infty} 
	\frac{\Pr( x_0 X_0 > x \mbox{ or } \ldots \mbox{ or } x_t X_t > x)}{\Pr(|X_0| > x)}
	} \\
	&=& \frac{1}{2}
	\E \bigl[ {\textstyle \bigvee_{j=0}^t \bigl\{x_j \sign(x_0) \lambda^{j/2} \prod_{i=1}^j Z_i \bigr\}_+^\alpha} \bigr],
\end{eqnarray*}
from which the domains of attraction of $(X_0, \ldots, X_t)$ in all $2^{t+1}$ corners of $\RR^{t+1}$ can be derived.
\end{example}

\begin{example}
\label{Ex:StochDiff}
Let $X_0, X_1, X_2, \ldots$ be defined recursively by the stochastic difference equation
\begin{equation}
\label{E:StochDiff}
	X_t = A_t X_{t-1} + B_t, \qquad t = 1, 2, \ldots
\end{equation}
where $(A_1, B_1), (A_2, B_2), \ldots$ are independent and identically distributed bivariate random vectors, independent of $X_0$. By \citet{Kesten73} and \citet[][Theorem~4.1]{Goldie91}, a stationary solution to \eqref{E:StochDiff} exists provided $\law( \log |A_1| \mid A_1 \neq 0)$ is non-lattice and there exists $0 < \alpha < \infty$ such that $\E [ |A_1|^\alpha ] = 1$, $\E [ |A_1|^\alpha (\log A_1)_+^\alpha ] < \infty$ and $\E [ |B_1|^\alpha ] < \infty$. Moreover, the stationary distribution satisfies
\[
	\lim_{x \to \infty} x^\alpha \Pr( \pm X_t > x) = C_\pm
\]
with $C_+$ and $C_-$ as in \citet[][equation~(4.3)]{Goldie91}, and $C_+ + C_- > 0$ if and only if $\Pr[ B_1 = (1 - A_1) c ] < 1$ for each $c \in \RR$. Obviously, if $A_1$ and $B_1$ are nonnegative, as in the case of a squared ARCH process \citep{dHRRdV89}, then $C_- = 0$ and thus $p = 1$.

Under these assumptions, Condition~\ref{C:RV} is satisfied with the given $\alpha$ and with $p = C_+ / (C_- + C_+)$, and Condition~\ref{C:phi} is satisfied with, in obvious notation, $\phi((a, b), \pm 1) = a$. The back-and-forth tail chain is given by Theorem~\ref{T:backforth}, the forward tail chain admitting the representation 
\[
	M_t = M_0 A_1 \cdots A_t, \qquad t = 1, 2, \ldots, 
\]
where $M_0, A_1, A_2, \ldots$ are independent random variables with $\Pr(M_0 = \pm 1) = C_\pm/(C_+ + C_-) = p(\pm 1)$ and $A_1, A_2, \ldots$ as in \eqref{E:StochDiff}; the backward tail chain, however, does in general not admit such a multiplicative random walk representation and its distribution is to be obtained through \eqref{E:BFTC:AB:4}.

By Theorem~\ref{T:MRV}, the joint distribution of $(X_0, \ldots, X_t)$ is multivariate regularly varying, and since $\E [ |M_1|^\alpha ] = \E [ |A_1|^\alpha ] = 1$, the limiting measure $\mu_t$ in \eqref{E:MRV} admits the simple expression in \eqref{E:simpler:b}. A useful special case is that for integer $t \geq 0$ and for non-zero $x_0, \ldots, x_t$,
\begin{eqnarray*}
	\lefteqn{
	\lim_{x \to \infty}
	\frac{\Pr( x_0 X_0 > x \mbox{ or } \ldots \mbox{ or } x_t X_t > x)}{\Pr(|X_0| > x)}
	} \\
	&=& p(\sign x_0)
	\E \bigl[ {\textstyle \bigvee_{j=0}^t \bigl\{ x_j \sign(x_0) \prod_{i=1}^j A_i \bigr\}_+^\alpha}\bigr].
\end{eqnarray*}
For nonnegative $A_t$ and positive $x_i$, this relation was already established in \citet[][Corollary~2.1]{GdHP04}; see also \citet{BDM02b}.
\end{example}

\begin{example}
\label{Ex:Smith}
Let $X_0, Z_1, J_1, Z_2, J_2, \ldots$ be independent random variables such that the common distribution function, $F$, of $X_0, Z_1, Z_2, \ldots$ is continuous and such that $J_1, J_2, \ldots$ are equal in distribution with $\Pr(J_t \in \{-1, 0, 1\}) = 1$. Denote $\Pr(J_t = j) = p_j$ for $j \in \{-1, 0, 1\}$ and assume $p_0 > 0$. Let $Q$ be the quantile function of $F$ and define $X_t$ for integer $t \geq 1$ recursively by
\begin{equation}
\label{E:Smith}
	X_t =
	\left\lbrace \begin{array}{l@{\quad}l}
		X_{t-1} & \mbox{if $J_t = 1$,} \\
		Z_t & \mbox{if $J_t = 0$,} \\
		Q(1-F(X_{t-1})) & \mbox{if $J_t = -1$.}
	\end{array} \right.
\end{equation}
The process $\{ X_t \}$ is a stationary Markov chain.

Now assume that the stationary marginal distribution $F$ satisfies Condition~\ref{C:RV} for some $0 < \alpha < \infty$ and some $0 < p < 1$. By the special form of the chain, Condition~\ref{C:phi} is then automatically fulfilled as well with
\[
	\phi((z,j),\pm 1) =
	\left\lbrace
		\begin{array}{l@{\quad}l}
		1 & \mbox{if $j = 1$,} \\
		0 & \mbox{if $j = 0$,} \\
		- \{(1-p)/p\}^{\pm 1/\alpha} &
		\mbox{if $j = -1$.}
		\end{array}
	\right.
\]
The back-and-forth tail chain of $\{ X_t \}$ is given by Theorem~\ref{T:backforth}; in terms of the representation in \eqref{E:BFTC:AB:1} to \eqref{E:BFTC:AB:4}, we have
\begin{eqnarray*}
	\law(A_{\pm 1}) &=& p_1 \delta_1 + p_0 \delta_0 + p_{-1} \delta_{-\{(1-p)/p\}^{1/\alpha}}, \\
	\law(B_{\pm 1}) &=& p_1 \delta_1 + p_0 \delta_0 + p_{-1} \delta_{-\{p/(1-p)\}^{1/\alpha}},
\end{eqnarray*}
see also Example~\ref{Ex:adjoint:Smith}. The fact that the forward and backward chain are equal in law is a consequence of the fact that the time-reversed process $\{ X_{-t} \}$ has the same distribution as the original one.
\end{example}

\begin{example}
\label{Ex:Smith:counterexample}
Consider again the chain $\{ X_t \}$ defined in \eqref{E:Smith} and assume that the stationary marginal distribution satisfies Condition~\ref{C:RV} for some $0 < \alpha < \infty$ but now with $p = 1$. Although the first part of Condition~\ref{C:phi}, equation~\eqref{E:phi:V}, is fulfilled with $\phi((z,j), +1) = 1$ if $j = 1$ and $\phi((z,j), +1) = 0$ if $j \in \{ -1, 0 \}$, the second part, equation~\eqref{E:phi:W}, is not, as $\lim_{x \to \infty} Q(1-F(-x)) / x = \infty$. The theory therefore breaks down: although
\[
	\law(X_1 / X_0 \mid X_0 > x) \dto \law(A_1) = p_1 \delta_1 + (p_0 + p_{-1}) \delta_0
\]
as $x \to \infty$, the forward tail chain is not equal to a multiplicative random walk with independent increments with common distribution $\law(A_1)$. Instead, if for instance $p_1 = 0$, then for integer $t \geq 1$ and as $x \to \infty$,
\[
	\law(X_1 / X_0, \ldots, X_{2t} / X_0 \mid X_0 > x) \dto
	\law \bigl(0, C_1, 0, C_1 C_2, \ldots, 0, {\textstyle \prod_{i=1}^t} C_i \bigr),
\]
where $C_1, C_2, \ldots$ are independent Bernoulli random variables with succes probability $p_{-1}^2$.
\end{example}

\begin{example}
\label{Ex:dependent}
%The following example is a counterexample to Theorem~\ref{T:forward} in case Condition~\ref{C:phi} is not fulfilled. If $p = 1$ in \eqref{E:p}, then the variables $M_0, M_1, M_2, \ldots$ in \eqref{E:Mj} and \eqref{E:forward:2} form a multiplicative random walk on $[0, \infty)$ given by $M_0 = 1$ and $M_t = \prod_{i=1}^t A_i$ where $A_1, A_2, \ldots$ are independent random variables with $\law(A_t) = \law(\phi(\eps_1, 1))$ for integer $t \geq 1$. However, if Condition~\ref{C:phi} is not fulfilled, then different kinds of limiting processes $\{ M_t \}$ may arise. 

Let $V_0, I_0, W_1, J_1, W_2, J_2, \ldots$ be independent random variables such that $\Pr(V_0 \leq x) = \Pr(W_t \leq x) = \exp(-1/x)$, $\Pr(I_0 = 1) = 1/2 = \Pr(I_0 = 0)$ and $\Pr(J_t = 1) = \pi = 1 - \Pr(J_t = 0)$ for real $x > 0$ and integer $t \geq 1$ and some $0 < \pi < 1$. Further, let $0 \leq a_0 < a_1 < 1$. Define $(V_t,I_t)$ for $t = 1, 2, \ldots$ recursively by
\begin{eqnarray*}
	V_t &=& \max \{ a_{I_{t-1}} V_{t-1}, (1 - a_{I_{t-1}}) W_t \}, \\
	I_t &=& 
	\left\lbrace \begin{array}{l@{\quad}l}
		I_{t-1} & \mbox{if $J_t = 1$,} \\
		1 - I_{t-1} & \mbox{if $J_t = 0$.}
	\end{array} \right.
\end{eqnarray*}
The process $\{ V_t \}$ may be thought of as a max-autoregressive process with parameter depending on a latent variable $I_t$. The process $\{ (V_t, I_t) \}$ forms a stationary Markov chain on $(0, \infty) \times \{0, 1\}$. Next, put
\[
	X_t = V_t + \floor{V_t} + I_t, \qquad t = 0, 1, 2, \ldots,
\]
where $\floor{v}$ is the integer part of a real $v$. The map $(v, i) \mapsto v + \floor{v} + i$ is a bijection from $(0, \infty) \times \{0, 1\}$ to $(0, \infty)$, the inverse map being $x \mapsto (v, i)$ where $v = \floor{x/2} + x - \floor{x}$ and $i$ is $0$ or $1$ depending on whether $\floor{x}$ is even or odd. Hence, $\{ X_t \}$ is a stationary Markov chain on $(0, \infty)$. 

Condition~\ref{C:RV} is satisfied with $\alpha = 1$ and $p = 1$. However, the first part of Condition~\ref{C:phi}, equation~\eqref{E:phi:V}, is not satisfied since the conditional distribution of $X_t$ given $X_{t-1}$ is strongly affected by whether $\floor{X_{t-1}}$ is odd or even through the value of $a_{I_{t-1}}$. The conclusion of Theorem~\ref{T:forward} is not satisfied either unless $\pi = 1/2$. In fact, for integer $t \geq 1$,
\[
	\law \biggl( \frac{X_0}{x}, \frac{X_1}{X_0}, \ldots, \frac{X_t}{X_0} \, \bigg| \, X_0 > x \biggr)
	\dto \law(Y, M_1, \ldots M_t)
\]
as $x \to \infty$, where $Y$ is independent of $M_1, M_2, \ldots$ and with distribution $\Pr(Y > y) = y^{-1}$ for $y \geq 1$, and where $M_t = A_1 \cdots A_t$ with $A_t = a_{I_t}$ for integer $t \geq 1$. Since the process $I_0, I_1, I_2, \ldots$ is a stationary Markov chain on $\{0, 1\}$ with transition matrix determined by $\Pr(I_t = I_{t+1}) = \pi$, the variables $A_1, A_2, \ldots$ form a stationary Markov chain on $\{ a_0, a_1 \}$. The variables $A_1, A_2, \ldots$ are not independent unless $\pi = 1/2$.
\end{example}

% ===============================================================================
\section*{Acknowledgement} 
The author thanks Bojan Basrak and Frank Redig for helpful discussions.

% ===============================================================================

\end{document}